\documentclass[copyright,creativecommons]{eptcs}
\providecommand{\event}{CCA 2010} 
\usepackage{breakurl}             

\usepackage{amsmath}
\usepackage{amsfonts}
\usepackage{graphicx}

\newtheorem{theorem}{Theorem}

\newtheorem{definition}{Definition}

\newtheorem{lemma}{Lemma}

\newtheorem{problem}{Problem}
\newtheorem{proposition}{Proposition}
\newtheorem{remark}{Remark}

\numberwithin{equation}{section}

\title{Computing
the speed of convergence of ergodic averages and pseudorandom points in
computable dynamical systems}

\author{Stefano Galatolo
\institute{Dipartiento di matematica applicata \\
Universita di Pisa}
\email{s.galatolo@ing.unipi.it}
\and
Mathieu Hoyrup  \qquad\qquad
\institute{LORIA, Vandoeuvre-l es-Nancy, France}
\email{Mathieu.Hoyrup@loria.fr}
\and
Crist\' obal Rojas \qquad\qquad
\institute{Fields Institute, Toronto, Canada}
\email{ cristobal.rojas@utoronto.ca}
}
\def\titlerunning{Speed of convergence of averages and pseudorandom points }
\def\authorrunning{S. Galatolo, M. Hoyrup, C Rojas }
\begin{document}
\maketitle

\begin{abstract}
A pseudorandom point in an ergodic dynamical system over a computable metric
space is a point which is computable but its dynamics has the same
statistical behavior as a typical point of the system.

It was proved in \cite{AvigadGT10} that in a system whose dynamics is
computable the ergodic averages of computable observables converge
effectively. We give an alternative, simpler proof of this result.

This implies that if also the invariant measure is computable then the
pseudorandom points are a set which is dense (hence nonempty) on the support
of the invariant measure.
\end{abstract}

\section{Introduction}
We will consider abstract algorithmic questions concerning the evolution of
a dynamical system. In particular, the algorithmic estimation of the speed
of convergence of ergodic averages and the recursive construction of points
whose dynamics is typical for the system.

This latter problem is related to the possibility of computer simulations,
as actual computers can only calculate the evolution of computable initial
conditions.

Let $X$ be a separable metric space and let the iterations of a map $T:X\rightarrow X$
define a dynamics. Let $\mu $ be an invariant measure for the dynamics ($\mu
(A)=\mu (T^{-1}(A))$ for each measurable set $A$). The famous Birkhoff
ergodic theorem says that if $\mu $ is ergodic (the only sets which are
invariant for the dynamics have full or null measure) then%
\begin{equation*}
\underset{n\rightarrow \infty }{\lim }\frac{1}{n}\sum f(T^{n}(x))=\int \!%
{f}\,\mathrm{d}{\mu }\text{ , }\mu -\text{almost everywhere.}
\end{equation*}

Similar results can be obtained for the convergence in the $L^{2}$ norm, and
others. The above result tells that, if the system is ergodic, there is a
fulll measure set of points for which the averages of the values of the
observable $f$ along its trajectory (time averages) coincides with the
spatial average of the observable $f$. Such a point could be called \textit{%
typical }for $f$ . Points which are typical for each $f$ which is continuous
with compact support are called typical for the measure $\mu $ (see
Definition \ref{mutyp}).

\noindent \textbf{Estimating the speed of convergence. }Many, more refined
results are linked to the speed of convergence of the above limit. We
consider this problem from the point of view of Computable Analysis. In the
paper \cite{AvigadGT10} some abstract results imply that in a computable
ergodic dynamical system, the speed of convergence of such averages can be
algorithmically estimated. On the other hand it is also shown that there are
non ergodic systems where this kind of estimations are not possible. In this
paper we present a simple, direct way to prove the following result:

\smallskip

\noindent \textbf{Theorem A. } \emph{If \ }$(X,\mu ,T)$\emph{\ is an ergodic
dynamical system and }$T$ \emph{is a.e. computable, then for each computable }$%
L^{1}$ \emph{observable }$f$\emph{\ the ergodic average }$A_{n}(x)=\frac{1}{n%
}\sum_{i=0}^{n-1}f(T^{i}(x))$ \emph{converge effectively a.e. to }$\int \!{f}%
\,\mathrm{d}{\mu }$ \emph{(see Definition \ref{def_rec_conv} for the precise
definition of this kind of convergence).} \newline

This theorem has some interesting consequences, as we are going to
illustrate.

\smallskip

\noindent \textbf{Computable points having typical statistical behavior}.
The set of computable points in $X$ (see Definition \ref{comp_points}) being
countable, is a very small (invariant) set, compared to the whole space. For
this reason, a computable point could be rarely be expected to be typical
for the dynamics, as defined before. More precisely, the Birkhoff ergodic
theorem and other theorems which hold for a full measure set, cannot help to
decide if there exist a computable point which is typical for the dynamics.
Nevertheless computable points are the only points we can use when we
perform a simulation or some explicit computation on a computer.

A number of theoretical questions arise naturally from all these facts. Due
to the importance of forecasting and simulation of a dynamical system's
behavior, these questions also have some practical motivation.

\begin{problem}
Since simulations can only start with computable initial conditions, given
some typical statistical behavior of a dynamical system, is there some
computable initial condition realizing this behavior? how to choose such
points?
\end{problem}

Such points could be called \emph{pseudorandom} points. Meaningful
simulations, showing typical behaviors of the dynamics can be performed if
computable, pseudorandom initial conditions exist ( and can be computed from the
description of the system).

We remark that it is widely believed that computer simulations produce
correct statistical behavior. The evidence is mostly heuristic. Most
arguments are based on the various \textquotedblleft
shadowing\textquotedblright\ results (see e.g. \cite{KH} chapter 18). In
this kind of approach (different from ours), it is possible to prove that in
a suitable system every pseudo-trajectory, as the ones which are obtained in
simulations with some computation error, is close to a real trajectory of
the system. However, even if we know that what we see in a simulation is
near to some real trajectory, we do not know if this real trajectory is
typical in some sense.

The main limit of this approach is that shadowing results hold only in
particular systems, having some uniform hyperbolicity, while many physically
interesting systems are not like this.

In our approach we consider real trajectories instead of "pseudo" ones and
we ask if there is some computable point which behaves as a typical point of
the space. Thanks to a kind of effective Borel-Cantelli lemma, in \cite%
{GHR07} the above problem is solved affirmatively for a class of systems
which satisfies certain technical assumptions which includes systems whose
decay of correlation is faster that $C\log ^{2}(time)$. In this paper we
prove the following more general result, as a consequence of the above
Theorem A:\newline

\noindent \textbf{Theorem B. }\emph{If \ }$(X,\mu ,T)$\emph{\ is a
computable dynamical system as above and }$\mu $\emph{\ is a computable
invariant ergodic measure,} \emph{there exist computable points }$\emph{x}$%
\emph{\ for which it holds:}%
\begin{equation}
\underset{n\rightarrow \infty }{\lim }\frac{f(x)+f(T(x))+\ldots
+f(T^{n-1}(x))}{n}=\int \!{f}\,\mathrm{d}{\mu }  \label{typic}
\end{equation}%
\emph{for any continuous function $f:X\rightarrow \mathbb{R}$ with compact
support.} \newline

The above theorem hence states that in such systems there are pseudorandom
points.

\noindent \textbf{Physical measures and computability. }To apply the above
corollary to concrete systems the main difficulty is to verify that the
invariant measure is computable. In \cite{GHR07} and \cite{GalHoyRoj3} it is
shown that this is verified for the physical\footnote{\textbf{\ }In general,
given $(X,T)$ there could be infinitely many invariant probability measures.
Among this class of measures, some of them are particularly important
because they are related to what can be seen in real experiments. Suppose
that we observe the behavior of the system $(X,T)$ through a class of
continuous functions $f_{i}:X\rightarrow \mathbb{R}$. We are interested in
the statistical behavior of $f_{i}$ along typical orbits of the system. Let
us suppose that the time average along the orbit of $x$ exists
\begin{equation*}
A_{x}(f_{i})=\underset{n\rightarrow \infty }{\lim }\frac{1}{n}\sum
f_{i}(T^{n}(x)).
\end{equation*}%
This is a real number for each $f_{i}.$ Moreover $A_{x}(f_{i})$ is linear
and continuous with respect to small changes of $f_{i}$ in the $\sup $ norm.
Then the orbit of $x$ acts as a measure $\mu _{x}$ and $A_{x}(f_{i})=\int \!{%
f_{i}}\,\mathrm{d}{\mu }_{x}$ (moreover this measure is also invariant for $%
T $). This measure is physically interesting if it is given by a
\textquotedblleft large\textquotedblright\ set of initial conditions. This
set will be called the basin of the measure. If $X$ is a manifold, it is
said that an invariant measure is \emph{physical} if its basin has positive
Lebesgue measure (see \cite{Y} for a survey and more precise definitions).}
invariant measure (the natural invariant measure to be considered in this
cases) in several classes of interesting systems as uniformly hyperbolic
systems, piecewise expanding maps and interval maps with an indifferent
fixed point. On the other hand there are computable systems having no
computable invariant measure ( \cite{GalHoyRoj3}) this shows some subtlety in this kind of
questions.

The way we handle computability on continuous spaces is largely inspired by
representation theory (see \cite{Wei00},\cite{BHW}). However, the main goal
of that theory is to study, in general topological spaces, the way
computability notions depend on the chosen representation. Since we focus
only on computable metric spaces, we do not use representation theory in its
general setting but instead present computability notions in a
self-contained way, and hopefully accessible to non-specialists.

\section{Computability}

The starting point of recursion theory was to give a mathematical definition
making precise the intuitive notions of \emph{algorithmic} or \emph{%
effective procedure} on symbolic objects. Several different formalizations
have been independently proposed (by Church, Kleene, Turing, Post,
Markov...) in the 30's, and have proved to be equivalent: they compute the
same functions from $\mathbb{N}$ to $\mathbb{N}$. This class of functions is
now called the class of \emph{recursive functions}. As an algorithm is
allowed to run forever on an input, these functions may be \emph{partial},
i.e. not defined everywhere. The \emph{domain} of a recursive function is
the set of inputs on which the algorithm eventually halts. A recursive
function whose domain is $\mathbb{N}$ is said to be \emph{total}.

We now recall an important concept from recursion theory. A set $E\subseteq
\mathbb{N}$ is said to be \textbf{\emph{recursively enumerable (r.e.)}} if
there is a (partial or total) recursive function $\varphi :\mathbb{N}%
\rightarrow \mathbb{N}$ enumerating $E$, that is $E=\{\varphi (n):n\in
\mathbb{N}\}$. If $E\neq \emptyset $, $\varphi $ can be effectively
converted into a total recursive function $\psi $ which enumerates the same
set $E$.

\subsection{Algorithms and uniform algorithms}

Strictly speaking, recursive functions only work on natural numbers, but
this can be extended to the objects (thought as ``finite'' objects) of any
countable set, once a numbering of its elements has been chosen. We will use
the word \emph{algorithm} instead of \emph{recursive function} when the
inputs or outputs are interpreted as finite objects. The operative power of
algorithms on the objects of such a numbered set obviously depends on what
can be effectively recovered from their numbers.

More precisely, let $X$ and $Y$ be numbered sets such that the numbering of $%
X$ is injective (it is then a bijection between $\mathbb{N}$ and $X$). Then
any \emph{recursive function} $\varphi:\mathbb{N}\to\mathbb{N} $ induces an
\emph{algorithm} $\mathcal{A}:X\to Y$. The particular case $X=\mathbb{N}$
will be much used.

For instance, the set $\mathbb{Q}$ of rational numbers can be injectively
numbered $\mathbb{Q}=\{q_0,q_1,\ldots\}$ in an \emph{effective} way: the
number $i$ of a rational $a/b$ can be computed from $a$ and $b$, and vice
versa. We fix such a numbering: from now and beyond the rational number with
number $i$ will be denoted by $q_i$.

Now, let us consider computability notions on the set $\mathbb{R}$ of real
numbers, introduced by Turing in \cite{Tur36}.

\begin{definition}
Let $x$ be a real number. We say that:

\begin{itemize}
\item $x$ is \textbf{\emph{lower semi-computable}} if the set $\{i\in
\mathbb{N}:q_i<x\}$ is r.e.

\item $x$ is \textbf{\emph{upper semi-computable}} if the set $\{i\in
\mathbb{N}:q_i>x\}$ is r.e.

\item $x$ is \textbf{\emph{computable}} if it is lower and upper
semi-computable.
\end{itemize}
\end{definition}

Equivalently, a real number is computable if and only if there exists an
algorithmic enumeration of a sequence of rational numbers converging
exponentially fast to $x$. That is:

\begin{proposition}
A real number is \textbf{\emph{computable}} if there is an algorithm $%
\mathcal{A}:\mathbb{N}\to \mathbb{Q}$ such that $|\mathcal{A}(n)-x|\leq
2^{-n}$ for all $n$.
\end{proposition}

\medskip

\noindent \textbf{Uniformity.} Algorithms can be used to define
computability notions on many classes of mathematical objects. The precise
definitions will be particular to each class of objects, but they will
always follow the following scheme:


\smallskip

\begin{center}
An object $O$ is \textbf{\emph{computable}} if there is an algorithm
\begin{equation*}
\mathcal{A}:X\to Y
\end{equation*}

which computes $O$ in some way.
\end{center}

\smallskip

Each computability notion comes with a uniform version. Let $(O_i)_{i\in%
\mathbb{N}}$ be a sequence of computable objects:

\smallskip

\begin{center}
$O_i$ is computable \textbf{\emph{uniformly in $\boldsymbol{i}$}} if there
is an algorithm
\begin{equation*}
\mathcal{A}:\mathbb{N}\times X \to Y
\end{equation*}

such that for all $i$, $\mathcal{A}_i:= \mathcal{A}(i,\cdot):X\to Y$
computes $O_i$.
\end{center}

\smallskip


For instance, the elements of a sequence of real numbers $(x_i)_{i\in\mathbb{%
N}}$ are uniformly computable if there is a algorithm $\mathcal{A}:\mathbb{N}%
\times\mathbb{N}\to\mathbb{Q}$ such that $|\mathcal{A}(i,n)-x_i|\leq 2^{-n}$
for all $i,n$.

In each particular case, the computability notion may take a particular
name: computable, recursive, effective, r.e., etc. so the term
``computable'' used above shall be replaced. 

\subsection{Computable metric spaces\label{CMS}}

A computable metric space is a metric space with an additional structure
allowing to interpret input and output of algorithms as points of the metric
space. This is done in the following way: there is a dense subset (called
ideal points) such that each point of the set is identified with a natural
number. The choice of this set is compatible with the metric, in the sense
that the distance between two such points is computable up to any precision
by an algorithm getting the names of the points as input. Using these simple
assumptions many constructions on metric spaces can be implemented by
algorithms.

\begin{definition}
A \textbf{\emph{computable metric space}} (CMS) is a triple $\mathcal{X}%
=(X,d,S)$, where

\begin{enumerate}
\item[(i)] $(X,d)$ is a separable metric space.

\item[(ii)] $S=\{s_{i}\}_{i\in \mathbb{N}}$ is a dense, numbered, subset of $%
X$ called the set of \textbf{\emph{ideal points}}.

\item[(iii)] The distances between ideal points $d(s_{i},s_{j})$ are all
computable, uniformly in $i,j$ (there is an algorithm $\mathcal{A}:\mathbb{N}%
^3\to \mathbb{Q}$ such that $|\mathcal{A}(i,j,n)-d(s_i,s_j)|<2^{-n}$).
\end{enumerate}
\end{definition}

$S$ is a numbered set, and the information that can be recovered from the
numbers of ideal points is their mutual distances. Without loss of
generality, we will suppose the numbering of $S$ to be injective: it can
always be made injective in an effective way.

We say that in a metric space $(X,d)$, a sequence of points $(x_{n})_{n\in
\mathbb{N}}$ converges\emph{\ recursively }to a point $x$ if there is an
algorithm $D:\mathbb{Q\rightarrow N}$ such that $d(x_{n},x)\leq \epsilon $
for all $n\geq D(\epsilon )$.

\begin{definition}
\label{comp_points}A point $x\in X$ is said to be \textbf{\emph{computable}}
if there is an algorithm $\mathcal{A}:\mathbb{N}\rightarrow S$ such that $(%
\mathcal{A}(n))_{n\in \mathbb{N}}$ converges recursively to $x$.
\end{definition}

We define the set of \textbf{\emph{ideal balls}} to be $\mathcal{B}%
:=\{B(s_{i},q_{j}):s_{i}\in S,0<q_{j}\in \mathbb{Q}\}$ where $B(x,r)=\{y\in
X:d(x,y)<r\}$ is an open ball. We fix a numbering $\mathcal{B}%
=\{B_{0},B_{1},\ldots \}$ which makes the number of a ball effectively
computable from its center and radius and vice versa. $\mathcal{B}$ is a
countable basis of the topology.

\begin{definition}[Effective open sets]
\label{reop} We say that an open set $U$ is \emph{effective} if there is an
algorithm $\mathcal{A}:\mathbb{N}\rightarrow \mathcal{B}$ such that $%
U=\bigcup_{n}\mathcal{A}(n)$.
\end{definition}

Observe that an algorithm which diverges on each input $n$ enumerates the
empty set, which is then an effective open set. Sequences of uniformly
effective open sets are naturally defined. Moreover, if $(U_{i})_{i\in
\mathbb{N}}$ is a sequence of uniformly effective open sets, then $%
\bigcup_{i}U_{i}$ is an effective open set.

\begin{definition}[Effective $G_{\protect\delta }$-set]
An \textbf{\emph{effective $\boldsymbol{G_{\delta }}$-set}} is an
intersection of a sequence of uniformly effective open sets.
\end{definition}

Obviously, an intersection of uniformly effective $G_{\delta }$-sets is also
an effective $G_{\delta }$-set.

Let $(X,S_{X}=\{s_{1}^{X},s_{2}^{X},...\},d_{X})$ and $(Y,S_{Y}=%
\{s_{1}^{Y},s_{2}^{Y},...\},d_{Y})$ be computable metric spaces. Let also $%
B_{i}^{X}$ and $B_{i}^{Y}$ be enumerations of the ideal balls in $X$ and $Y$%
. A computable function $X\rightarrow Y$ is a function whose behavior can be
computed by an algorithm up to any precision. For this it is sufficient that
the pre-image of each ideal ball can be effectively enumerated by an
algorithm.

\begin{definition}[Computable Functions]
\label{comp_func} A function $T:X\rightarrow Y$ is \textbf{\emph{computable}}
if $T^{-1}(B_{i}^{Y})$ is an effective open set, uniformly in $i$. That is,
there is an algorithm $\mathcal{A}:\mathbb{N}\times \mathbb{N}\rightarrow
\mathcal{B}^{X}$ such that $T^{-1}(B_{i}^{Y})=\bigcup_{n}\mathcal{A}(i,n)$
for all $i$.

\noindent A function $T:X\rightarrow Y$ is \textbf{\emph{computable on $%
D\subseteq X$}} if there are uniformly effective open sets $U_{i}$ such that
$T^{-1}(B_{i}^{Y})\cap D=U_{i}\cap D.$
\end{definition}

\subsection{Computable measures\label{seccompmu}}

Let us consider the space $PM(X)$ of Borel probability measures over $X$.
Let $C_{0}(X)$ be the set of real-valued bounded continuous functions on $X$%
. We recall the notion of weak convergence of measures:

\begin{definition}
$\mu _{n}$ is said to be \emph{\textbf{weakly convergent}} to $\mu $ if $%
\int\!{f}\,\mathrm{d}{\mu_n}\rightarrow \int\!{f}\,\mathrm{d}{\mu}$ for each
$f\in C_{0}(X)$.
\end{definition}

Let us introduce the Wasserstein-Kantorovich distance between measures. Let $%
\mu _{1}$ and $\mu _{2}$ be two probability measures on $X$ and consider:

\begin{equation*}
W_{1}(\mu _{1},\mu _{2})=\underset{f\in 1\text{-Lip}(X)}{\sup }\left|\int\!{f%
}\,\mathrm{d}{\mu_1}-\int\!{f}\,\mathrm{d}{\mu _2}\right|,
\end{equation*}

\noindent where $1\mbox{-Lip}(X)$ is the space of functions on $X$ having
Lipschitz constant less than one.

\begin{proposition}[see \protect\cite{AGS} Prop 7.1.5]
\label{ambros}\mbox{}

\begin{enumerate}
\item $W_{1}$ is a distance and if $X$ is bounded, separable and complete,
then $PM(X) $ with this distance is a separable and complete metric space.

\item If $X$ is bounded, a sequence is convergent for the $W_{1}$ metrics if
and only if it is convergent for the weak topology.
\end{enumerate}
\end{proposition}

Item (1) has an effective version: $PM(X)$ inherits the computable metric
structure of $X$. Indeed, given the set $S _{X}$ of ideal points of $X$ we
can naturally define a set of ideal points $S _{PM(X)}$ in $PM(X)$ by
considering finite rational convex combinations of the Dirac measures $%
\delta _{s}$ supported on ideal points $s\in S_{X}$. This is a dense subset
of $PM(X)$. The proof of the following proposition can be found in (\cite%
{HR07}).

\begin{proposition}
If $X$ bounded then $(PM(X),W_1, S _{PM(X)})$ is a computable metric space.
\end{proposition}

A measure $\mu $ is then computable if there is a sequence $\mu _{n}\in S_{PM(X)}$ converging exponentially fast to $\mu $ in the $W_{1}$ metric (and
hence $\mu _{n}$ weakly converge to $\mu $).

\subsection{Computable probability spaces}

To obtain computability results on dynamical systems, it seems obvious that
some computability conditions must be required on the system. The
\textquotedblleft good\textquotedblright\ conditions, if any, are not
obvious to specify.

A computable function defined on the whole space is necessarily continuous.
But a transformation or an observable need not be continuous at every point,
as many interesting examples prove (piecewise-defined transformations,
characteristic functions of measurable sets,... ), so the requirement of
being computable everywhere is too strong. In a measure-theoretical setting,
a natural weaker condition is to require the function to be computable on a
set of full measure. It can be proved that such a function can be extended
to a function which is computable on a full-measure effective $G_{\delta }$%
-set (see \cite{HR07, Hoy08}).

\begin{definition}
\label{a.e.compufunct}A \textbf{\emph{computable probability space}} is a
pair $(X,\mu )$ where $X$ is a computable metric space and $\mu $ a
computable Borel probability measure on $X$.

Let $Y$ be a computable metric space. A function $(X,\mu )\rightarrow Y$ is
\textbf{\emph{almost everywhere computable}} (a.e. computable for short) if
it is computable on an effective $G_{\delta }$-set of measure one, denoted
by $\mathrm{dom}f$ and called the \emph{domain of computability of $f$}.

A \textbf{\emph{morphism}} of computable probability spaces $%
f:(X,\mu)\to(Y,\nu)$ is a morphism of probability spaces which is a.e.
computable.
\end{definition}

\begin{remark}
\label{com_sum} A sequence of functions $f_{n}$ is uniformly a.e. computable
if the functions are uniformly computable on their respective domains, which
are uniformly effective $G_{\delta }$-sets. Observe that in this case,
intersecting all the domains provides an effective $G_{\delta }$-set on
which all $f_{n}$ are computable. In the following we will apply this
principle to the iterates $f_{n}=T^{n}$ of an a.e. computable function $%
T:X\rightarrow X$, which are uniformly a.e. computable.
\end{remark}

\begin{remark}\label{l1comp}
The space $L^1(X,\mu)$ (resp. $L^2(X,\mu)$) can be made a computable metric
space, choosing some dense set of bounded computable functions as ideal
elements. We say that an integrable function $f:X\to\overline{\mathbb{R}}$
is $L^1(X,\mu)$-computable if its equivalence class is a computable element
of the computable metric space $L^1(X,\mu)$. Of course, if $f=g$ $\mu$-a.e.,
then $f$ is $L^1(X,\mu)$-computable if and only if $g$ is. Basic operations
on $L^1(X,\mu)$, such as addition, multiplication by a scalar, $min$, $max$ etc. are computable. Moreover, if $T:X\to X$ preserves $\mu$ and $T$ is a.e. computable, then $f\to f\circ T$ (from $L^1$ to $L^1$) is computable  (see \cite{HoyRojCiE09}).
\end{remark}

\subsubsection{Application to convergence of random variables}

Here, $(X,\mu)$ is a computable probability space, where $X$ is complete.

\begin{definition}
A \textbf{\emph{random variable}} on $(X,\mu)$ is a measurable function $%
f:X\to \mathbb{R}$.
\end{definition}

\begin{definition}
Random variables $f_n$ \textbf{\emph{effectively converge in probability}}
to $f$ if for each $\epsilon>0$, $\mu\{x: |f_n(x)-f(x)|<\epsilon\}$
converges effectively to $1$, uniformly in $\epsilon$. That is, there is a
computable function $n(\epsilon,\delta)$ such that for all $n\geq
n(\epsilon,\delta)$, $\mu\{|f_n-f|\geq\epsilon\}<\delta$.
\end{definition}

\begin{definition}
\label{def_rec_conv}Random variables $f_{n}$ \textbf{\emph{effectively
converge almost surely}} to $f$ if $f_{n}^{\prime }=\sup_{k\geq n}|f_{k}-f|$
effectively converge in probability to $0$.
\end{definition}

The following result (\cite{GHR07}, Theorem 2) shows that if a sequence $%
f_{n}$ converges effectively a.s. to $f$ then there are computable points
which for which $f_{n}(x)\rightarrow f(x)$.

\begin{theorem}
\label{theorem_convergence_Borel_Cantelli} Let $X$ be a complete metric
space. Let $f_{n},f$ be uniformly a.e. computable random variables. If $%
f_{n} $ effectively converges almost surely to $f$ then the set $%
\{x:f_{n}(x)\rightarrow f(x)\}$ contains an effective Borel-Cantelli set
(see the Appendix for the precise definition).

In particular, it contains a sequence of uniformly computable points which
is dense in $\mathrm{Supp}(\mu )$.
\end{theorem}

\begin{remark}
\label{remarkmistery}Moreover, the effective Borel Cantelli Set found above
depends algorithmically on $f_{n}$ and on the function $n(\delta,\epsilon)$
giving the rate of convergence (see the proof of Theorem 2 in \cite{GHR07}).
Hence the result is uniform in $f_n$ and $n(\delta,\epsilon)$.
\end{remark}

\subsection{Effective $L^{1}$, $L^{2}$ convergence}

Let $(X,\mu ,T)$ be a computable measure-preserving system and $f$ a $L^{1}$%
-computable function (in the sense that it is a computable point of the
metric space $L^{1}$). It was proved in \cite{AvigadGT10} that the ($
L^{1},L^{2}$ and almost sure) convergence of the Birkhoff averages of $f$ is
effective as soon as the norm of the limit $f^{\ast }$ is computable. Here we give an
alternative proof in the ergodic case which is simpler as it uses the classical convergence
result instead of giving a \textquotedblleft constructive\textquotedblright\
proof.

Let us call $(X,\mu ,T)$ a \textbf{\emph{computable ergodic system}} if $%
(X,\mu )$ is a computable probability space where $T$ is an endomorphism
(i.e. an a.e. computable measure-preserving transformation) and $(X,\mu ,T)$
is ergodic. Let $||f||$ denote the $L^{1}$ norm or the $L^{2}$ norm.

\begin{proposition}
Let $(X,\mu ,T)$ be a computable ergodic system. Let $f$ be a computable
element of $L^{1}(X,\mu )$ (resp. $L^{2}(X,\mu )$).

The $L^1$ convergence (resp. $L^2$ convergence) of the Birkhoff averages of $%
f$ is effective.
\end{proposition}

{\bf Proof.}
Replacing $f$ with $f-\int fd\mu $, we can assume that $\int fd\mu =0$.
Let $A_{n}=(f+f\circ T+\ldots +f\circ T^{n-1})/n$. The sequence $||A_{n}||$
is computable ( see Remark \ref{l1comp} ) and converges to $0$ by the ergodic
theorems.

Given $p\in N$, we write $m\in N$ as $m=np+k$ with $0\leq k<p$. Then

\begin{eqnarray*}
A_{np+k} &=&\frac{1}{np+k}\left( \sum_{i=0}^{n-1}pA_{p}\circ
T^{pi}+kA_{k}\circ T^{pn}\right) \\
||A_{np+k}|| &\leq &\frac{1}{np+k}(np||A_{p}||+k||A_{k}||) \\
&\leq &||A_{p}||+\frac{||{A_{k}}||}{n} \\
&\leq &||{A_{p}}||+\frac{||{f}||}{n}.
\end{eqnarray*}

Let $\epsilon >0$. We can compute some $p=p(\epsilon )$ such that $||{A_{p}}%
||<\epsilon /2$. Then we can compute some $n(\epsilon )\geq \frac{2}{%
\epsilon }||{f}||$. The function $m(\epsilon ):=n(\epsilon )p(\epsilon )$ is
computable and for all $m\geq m(\epsilon )$, $||{A_{m}}||\leq \epsilon $.
$\Box$

\subsection{Effective almost sure convergence}

Now we use the above result to find a computable estimation for the a.s.
speed of convergence.

\begin{theorem}
\label{corollaryx}Let $(X,\mu ,T)$ be a computable ergodic system. If $f$ is
$L^{1}(X,\mu )$-computable, then the \ a.s. convergence is effective.
\end{theorem}

This will be proved by the following

\begin{proposition}
\label{theorem_ergodic_effective_as_convergence} If $f$ is $L^{1}(X,\mu )$%
-computable, and $||{f}||_{\infty }$ is bounded, then
the almost-sure convergence is effective (uniformly in $f$ and a bound on  $||{f}||_{\infty }$). 
\end{proposition}

To prove this we will use the Maximal ergodic theorem:

\begin{lemma}[Maximal ergodic theorem]
\label{lemma_maximal} For $f\in L^{1}(X,\mu )$ and $\delta >0$,
\begin{equation*}
\mu (\{\sup_{n}|A_{n}^{f}|>\delta \})\leq \frac{1}{\delta }||{f}||_{1}.
\end{equation*}
\end{lemma}

The idea is simple: compute some $p$ such that $||A_p^f||_{1}$ is small,
apply the maximal ergodic theorem to $g:=A_{p}^{f}$, and then there is $%
n_{0} $, that can be computed, such that $A_{n}^{f}$ is close to $A_{n}^{g}$
for $n\geq n_{0}$.

{\bf Proof.}
Let $\epsilon ,\delta >0$. Compute $p$ such that $||{A_{p}^{f}}||\leq \delta
\epsilon /2$. Applying the maximal ergodic theorem to $g:=A_{p}^{f}$ gives:
\begin{equation}
\mu (\{\sup_{n}|A_{n}^{g}|>\delta /2\})\leq \epsilon .  \label{bound_g}
\end{equation}

Now, $A_{n}^{g}$ is not far from $A_{n}^{f}$: expanding $A_{n}^{g}$, one can
check that
\begin{equation*}
A_{n}^{g}=A_{n}^{f}+\frac{u\circ T^{n}-u}{np},
\end{equation*}%
where $u=(p-1)f+(p-2)f\circ T+\ldots +f\circ T^{p-2}$. $||{u}||_{\infty
}\leq \frac{p(p-1)}{2}||{f}||_{\infty }$ so if $n\geq n_{0}\geq 4(p-1)||{f}%
||_{\infty }/\delta $, then $||{A_{n}^{g}-A_{n}^{f}}||_{\infty }\leq \delta
/2$. As a result, if $|A_{n}^{f}(x)|>\delta $ for some $n\geq n_{0}$, then $%
|A_{n}^{g}(x)|>\delta /2$. From (\ref{bound_g}), we then derive
\begin{equation*}
\mu (\{\sup_{n\geq n_{0}}|A_{n}^{f}|>\delta \})\leq \epsilon .
\end{equation*}

As $n_{0}$ can be computed from $\delta $ and $\epsilon $, we get the result.
$\Box $

\begin{remark}
\label{remarkboundunif}This result applies uniformly to a uniform sequence
of computable $L^{\infty }(X,\mu )$ observables $f_{n}$.
\end{remark}

We now extend this to $L^{1}(X,\mu )$-computable functions, using the
density of $L^{\infty }(X,\mu )$ in $L^{1}(X,\mu )$.

{\bf Proof.}
(of Theorem \ref{corollaryx}) Let $\epsilon ,\delta >0$. For $M\in \mathbb{N}
$, let us consider $f_{M}^{\prime }\in L^{\infty }(X,\mu )$ defined as
\begin{equation*}
f_{M}^{\prime }(x)=\left\{
\begin{array}{cc}
\min (f,M) & if~f(x)\geq 0 \\
\max (f,-M) & if~f(x)\leq 0.%
\end{array}%
\right.
\end{equation*}

Compute $M$ such that $||{f-f_{M}^{\prime }}||_{1}\leq \delta \epsilon $.
Applying Proposition \ref{theorem_ergodic_effective_as_convergence} to $%
f_{M}^{\prime }$ gives some $n_{0}$ such that $\mu (\{\sup_{n\geq
n_{0}}|A_{n}^{f_{M}^{\prime }}|>\delta \})<\epsilon $. Applying Lemma \ref%
{lemma_maximal} to $f_{M}^{\prime \prime }=f-f_{M}^{\prime }$ gives $\mu
(\{\sup_{n}|A_{n}^{f_{M}^{\prime \prime }}|>\delta \})<\epsilon $. As a
result, $\mu (\{\sup_{n\geq n_{0}}|A_{n}^{f}|>2\delta \})<2\epsilon $.
$\Box $
\begin{remark}
\label{remarkboundunif2} Also Theorem \ref{corollaryx} applies uniformly on
an uniform sequence of computable $L^{1}(X,\mu )$ observables $f_{n}$.
\end{remark}

\begin{remark}
\label{remarkcompuL1} We remark that a bounded a.e. computable function, as
defined in Definition \ref{a.e.compufunct} is a computable element of $%
L^{1}(X,\mu )$ (see \cite{HoyRojCiE09}). Conversely, if $f$ is a computable
element of $L^{1}(X,\mu )$ then there is a sequence of uniformly computable
functions $f_{n}$ that effectively converge $\mu $-a.e. to $f$.
\end{remark}

\section{Pseudorandom points and dynamical systems}

As said before the famous Birkhoff ergodic theorem says that in an ergodic
system, the time average computed along $\mu $-almost every orbit coincides
with space average with respect to $\mu $.

If a point $x$ satisfies Equation \ref{typic} for a certain $f$, then we say
that $x$ is typical with respect to the observable $f$.

\begin{definition}
\label{mutyp}A point $x$ is $\boldsymbol{\mu}$\textbf{\emph{-typical}} if $x$
is typical w.r.t. every continuous function $f:X\to \mathbb{R}$ with compact
support.
\end{definition}

We remark that \emph{from now on we will suppose that }$X$\emph{\ is a
complete metric space}. We will see that such $\mu -$typical points exist in
computable ergodic systems. First we give a result for $L^{1}$ observables.

\begin{theorem}
\label{thm1}If $(X,\mu ,T)$ is a computable ergodic system, $f$ is $%
L^{1}(X,\mu )$ and a.e. computable then there is a uniform sequence $x_{i}$
of computable points which is dense on the support of $\mu $ such that for
each $i$
\begin{equation*}
\underset{n\rightarrow \infty }{\lim }\frac{1}{n}\sum f(T^{n}(x_{i}))=\int \!%
{f}\,\mathrm{d}{\mu .}
\end{equation*}
\end{theorem}

{\bf Proof.}
Apply theorem \ref{theorem_convergence_Borel_Cantelli} to the sequence of
uniformly a.e. computable functions $f_{n}={A_n}^f$ which
converge effectively almost-surely by theorem \ref{corollaryx}. We obtain
that the set of points for which $\frac{1}{n}\sum f(T^{n}(x_{n}))\rightarrow
\int \!{f}\,\mathrm{d}{\mu }$ contains a sequence of computable points, as
in the statement.
$\Box $

Let $g$ be the point-wise limit of a sequence of uniformly computable
functions $f_{n}$ (defined $\mu $-a.e., as in Remark \ref{remarkboundunif2}%
): Theorem \ref{thm1} can also be proved to hold for the observable $g$,
i.e. in a computable ergodic system there exists computable points for which
the Birkhoff averages of $g$ converge to $\int \!{g}\,\mathrm{d}{\mu }$.

Since it is possible to construct a r.e. set of computable functions which
is dense in the space of compactly supported continuous functions we can
also obtain the following

\begin{theorem}
If $(X,\mu ,T)$ is a computable ergodic system then there is a uniform
sequence $x_{n}$ of computable points which in dense on the support of $\mu $
such that for each $n$, $x_{n}$ is $\mu -$typical.
\end{theorem}

{\bf Proof.}
Let us introduce (following \cite{Gac05}) a certain fixed, enumerated
sequence of Lipschitz functions. Let $\mathcal{F}_{0}$ be the set of
functions of the form:
\begin{equation}
g_{s,r,\epsilon }=|1-|d(x,s)-r|^{+}/\epsilon |^{+}  \label{lip functions}
\end{equation}%
where $s\in S$, $r,\epsilon \in \mathbb{Q}$ and $|a|^{+}=\max \{a,0\}$.

$g_{s,r,\epsilon }$ is a Lipschitz functions whose value is 1 inside the
ideal ball $B(s,r)$, 0 outside $B(s,r+\epsilon )$ and with intermediate
values in between. It is easy to see that the real-valued functions $%
g_{s_{i},r_{j},\epsilon _{k}}:X\rightarrow \mathbb{R}$ are computable,
uniformly in $i,j,k$.

Let $\mathcal{F}$ be the smallest set of functions containing $\mathcal{F}%
_{0}$ and the constant 1, and closed under $\max $, $\min $ and rational
linear combinations. Clearly, this is also a uniform family of computable
functions. We fix some enumeration $\nu _{\mathcal{F}}$ of $\mathcal{F}$ and
we write $g_{n}$ for $\nu _{\mathcal{F}}(n)\in \mathcal{F}$. We remark that
this set is dense in the set of continuous functions with compact support.

By Remark \ref{remarkcompuL1} moreover these are computable elements of $%
L^{1}(X,\mu )$, hence Theorem \ref{corollaryx} applies uniformly to these
observables. This means that we can apply Theorem \ref%
{theorem_convergence_Borel_Cantelli} uniformly on this sequence. By
intersecting all effective Borel-Cantelli sets given by Theorem \ref%
{theorem_convergence_Borel_Cantelli}, since the intersection of a uniform
family of effective BC sets contains an effective BC set (see Remark \ref%
{remarkmistery} and Proposition \ref{intersection_BC_proposition}) and such
a set contains a sequence of computable points which are dense in the
support of $\mu $ (see Theorem \ref{effective_BC_theorem}), which are
typical {\em for all} observables in $\mathcal{F}$ (in the same way as in Theorem %
\ref{thm1}). Since each continuous function $f$ with compact support can be
approximated in the $L^{\infty }$ norm by a function in $\mathcal{F}$ (and
in particular the approximating function in $\mathcal{F}$ have values near
the values of $f$ at almost each point) the statement is proved.
$\Box $
\subsection{Conclusion and some open question}

We have seen that in computable ergodic systems, the speed of a.e.
convergence of ergodic averages is computable, and if moreover also the
invariant measure is computable then there are computable points which are
typical for the statistical behavior of the system: the pseudorandom points.

The assumption about the computability of the measure is not redundant with
the assumption about the computability of the system because, as said
before, there are computable systems having not computable invariant
measures. It is also interesting to remark that there 
are systems where the map $T$ is computable, for which we can consider a non computable ergodic measure having no pseudorandom points (an example can be constructed considering a Bernoulli shift where the symbols have non computable probability).

Our results about the existence of these points, however, tell not much
about the computational complexity which is necessary to find them. It would
be of practical importance to have fast algoritms for this computation.

All these questions are related to another general (and vague) question we
like to cite, which is of great practical importance:  why, many "naive"
simulations of  dynamical systems give reasonable results? or more
precisely: under which assumptions the simulations of a dynamical system by
a computer (which is a kind of discrete model for a continuous phenomena)
give correct results?

\section{Appendix: effective BC sets}

We recall some results from \cite{GHR07} which are used in the proofs of the
present paper.

 Given a measurable space $X$ endowed with a probability
measure $\mu $, the well known Borel-Cantelli lemma states that if a
sequence of sets $A_{k}$ is such that $\sum \mu (A_{k})<\infty $ then the
set of points which belong to finitely many $A_{k}$'s has full measure. It holds that if the $A_{k}$ are given in some \textquotedblleft
effective\textquotedblright\ way (and $\mu $ is computable) then this full
measure set contains computable points, which can be effectively constructed.

\begin{definition}
A sequence of positive numbers $a_{i}$ is \emph{effectively summable} if the
sequence of partial sums converges effectively: there is an algorithm $%
\mathcal{A}:\mathbb{Q\rightarrow N}$ such that if $A(\epsilon )=n$ then $%
\sum_{i\geq n}a_{i}\leq \epsilon $.
\end{definition}

For the sake of simplicity, we will focus on the complements $U_n$ of the $%
A_n$.

\begin{definition}
\label{ebc}An \emph{effective Borel-Cantelli sequence} is a sequence $%
(U_{n})_{n\in \mathbb{N}}$ of uniformly effective open sets such that the
sequence $\mu(X\setminus U_{n})$ is effectively summable.

The corresponding \emph{effective Borel-Cantelli set} is $%
\bigcup_{k}\bigcap_{n\geq k}U_{n}$.
\end{definition}

\begin{proposition}
\label{intersection_BC_proposition} The intersection of any uniform family
of effective Borel-Cantelli sets contains an effective Borel-Cantelli set.
\end{proposition}

\begin{theorem}
\label{effective_BC_theorem} Let $X$ be a complete Computable Metric Space
and $\mu $ a computable Borel probability measure on $X$.

Every effective Borel-Cantelli set $R$, contains a sequence of uniformly
computable points which is dense in the support of $\mu$.
\end{theorem}

\end{document}